\documentclass[12pt]{amsart}

\usepackage[T1]{fontenc}        
\usepackage{ae,aecompl}
\usepackage{enumerate}
\usepackage{latexsym}
\usepackage{amsmath,amsthm,amsfonts}
\usepackage{amssymb}
\usepackage{epic}
\usepackage{tikz}
\usetikzlibrary{calc}
\usepackage{graphicx}
\usepackage{hyperref} 
\usepackage{url}
\usepackage{pgfplots}
\usepackage{subcaption}
\usepackage{algorithm}
\usepackage[noend]{algpseudocode}
\usepackage{float}
\usepackage{multicol}
\usetikzlibrary{matrix,shapes,decorations.pathreplacing}
\usepackage{easybmat}
\usepackage{multirow,bigdelim}
\makeatletter
\@namedef{subjclassname@2020}{%
  \textup{2020} Mathematics Subject Classification}
\makeatother

\usepackage[backend=biber]{biblatex} 
\DeclareFieldFormat{pages}{#1}
\title[Progressions in Euclidean Ramsey theory]{Progressions in Euclidean Ramsey theory}

\author[]{Jakob F\"uhrer}
\address{Jakob F\"uhrer, Institute of Analysis and Number Theory, Graz University of Technology,
Kopernikusgasse 24/II,
8010 Graz, Austria.}
\email{jakob.fuehrer@tugraz.at}

\author[]{G\'eza T\'oth}
\address{G\'eza T\'oth, R\'enyi A. Institute of Mathematics, Re\'altanoda u. 13-15, H-1053, Budapest, Hungary}
 \email{geza@renyi.hu}
\date{}

\usepackage{geometry}

\usepackage{tikz-3dplot}
\usetikzlibrary{calc}
\pgfplotsset{compat=1.18} 

\begin{document}
\baselineskip=17pt

 \maketitle
\newtheorem{theorem}{Theorem}
\newtheorem*{klar}{Klar}
\newtheorem{method}{Method}

\newtheorem{lemma}{Lemma}

\newtheorem{cor}{Corollary}

\newtheorem{conjecture}{Conjecture}

\theoremstyle{definition}
\newtheorem{defi}{Definition}

\newtheorem{bsp}{Beispiel}

\newtheorem*{bem}{Bemerkung}

\newtheorem*{vorschau}{Vorschau}

\newtheorem*{erg}{Ergänzung}

\theoremstyle{remark}
\newtheorem{remark}{Remark}

\newtheorem*{notation}{Notation}

\newtheorem{claim}{Claim}[theorem]
\renewcommand{\theclaim}{\arabic{claim}}
\newenvironment{proofofclaim}[1][\proofname\ of Claim \theclaim]{%
  \proof[#1]%
  \renewcommand\qedsymbol{$\blacksquare$}
}{\endproof}

\newcommand{\ndiv}{\not \hspace{3pt} \mid }

\newcommand*\hexbrace[2]{%
  \underset{#2}{\underbrace{\rule{#1}{0pt}}}}

\section*{Abstract}
\label{cha:abstract}

Conlon and Wu \cite{conlon2022lines} showed that there is a red/blue-coloring of $\mathbb{E}^n$ that does not contain $3$ red collinear points separated by unit distance
and $m=10^{50}$ blue collinear points separated by unit distance. 
We prove that the statement holds with $m=1177$. We show similar  
results with different distances between the points.

\section{Introduction}

For $n\in\mathbb{N}$ let $\mathbb{E}^n$ the $n$-dimensional Euclidean space, i.e. $\mathbb{R}^n$ with Euclidean distances. For any $m>0$ integer, let $\ell_m$ denote the {\em $m$-progression} with distance $1$, that is, a set of $m$ points on a line 
so that there is a unit distance between the consecutive points. 
In general, for any $\alpha\in\mathbb{R}_+$, $\alpha\ell_m$ stands for 
an $m$-progression with distance $\alpha$, 
that is, a set of $m$ points on a line 
so that there is a distance $\alpha$ between the consecutive points. 

For any finite sets $A, B\subset \mathbb{E}^n$, 
we write $\mathbb{E}^n\rightarrow (A, B)$ if for every red/blue-coloring of $\mathbb{E}^n$, there is either a red copy of $A$ or a blue copy of $B$.
Conversely,  write $\mathbb{E}^n\nrightarrow (A, B)$ if a red/blue-coloring of $\mathbb{E}^n$ exists that does not contain any red copy of $A$ nor any blue copy of $B$.

In this note we investigate the case where $A$ and $B$ are progressions. The general question is that 
for which $n$, $m_1$, $m_2$ does $\mathbb{E}^n\rightarrow (\ell_{m_1},\ell_{m_2})$ hold.
These kind of problems, in a much more general form, were first studied in a series of papers by 
Erd\H os, Graham, Montgomery, Rothschild, 
Spencer and Straus \cite{EGM1, EGM2, EGM3}.

Conlon and Fox \cite{ConlonFox} proved that there is a constant $c>0$ such that $\mathbb{E}^n\nrightarrow (\ell_2,\ell_m)$ for all $m\geq 2^{cn}$. 
However, it follows from a result of Szlam \cite{S01} and Frankl and Wilson \cite{FW81} that $\mathbb{E}^n\rightarrow (\ell_2,\ell_m)$ for some other constant $c'>0$ and  $m\leq 2^{c'n}$.

If we replace $\ell_2$ with $\ell_3$, the situation is quite different. 
Conlon and Wu \cite{conlon2022lines} showed that $\mathbb{E}^n\nrightarrow (\ell_{3},\ell_{m})$ for $m=10^{50}$, independent of the dimension $n$.
The proof is based on a  random spherical coloring, that is, 
the color of each point depends only on its distance from the origin. So it can be applied in any dimension.
Spherical colorings were applied by  Erd\H os et al.  \cite{EGM1} to show a four-coloring of $E^n$ with no monochromatic $\ell_3$. In this note we construct an explicit spherical coloring to improve the bound 
$10^{50}$ to $1177$.


\begin{theorem}
\label{mainthm}
For any $n>0$, there exists a red/blue-coloring of $\mathbb{E}^n$ that does not contain any red copy of $\ell_3$ and any blue copy of $\ell_{1177}$.
\end{theorem}

We also studied progressions with different distances. 
Observe that the statement $\mathbb{E}^n\rightarrow (\alpha_{red} \ell_{m_1},\alpha_{blue} \ell_{m_2})$ is equivalent to $\mathbb{E}^n\rightarrow ( \ell_{m_1}, (\alpha_{blue}/\alpha_{red})\ell_{m_2})$.

\begin{theorem}
\label{generalthm}

For any $n>0$, there exists a red/blue-coloring of $\mathbb{E}^n$ that does not contain any red copy of $\ell_3$ and any blue copy of $\alpha\ell_{8649}$, whenever $\alpha\in\mathbb{R}_+$ satisfies at least one of the following conditions:
\begin{itemize}
\item $\alpha^2\not\in \mathbb{Q}$,
\item $\alpha^2=p/q$, $p,q\in\mathbb{N}$ and $47\not\vert \; q$,
\item $\alpha^2\geq 2$,
\item $\alpha^2\leq 1/(7\cdot 47^4\cdot 48).$
\end{itemize}
\end{theorem}


For integers $a\le b$ let $[a,b]=\{a,a+1,...,b\}$. 
For any real $\gamma$ let $\lfloor \gamma \rfloor$ be the integral part of $\gamma$ and let $\{\gamma\}:=\gamma-\lfloor \gamma \rfloor$ be the fractional part.
Denote by $\mathbb{F}_p$ the field with $p$ elements and by $S_p$ the set of squares in $\mathbb{F}_p$. 
We write $\mathbb{F}_p^*$ and $S_p^*$ for the corresponding sets without zero.

\section{Overview}

In Section \ref{mainproof} and \ref{generalproof} we prove Theorems~\ref{mainthm} and \ref{generalthm}, respectively. 
Both proofs follow the same ideas which we sketch in this section.

We start with the most important fact: If three points $x,y,z$ lie in arithmetic progression, their norms satisfy an equation of the form 
\begin{equation*}
\label{lineq}
|x|^2-2|y|^2+|z|^2=K,  
\end{equation*}
where $K=2|x-y|^2=2|y-z|^2$, so it depends on the distance of the points and  
not on their locations.  
It is therefore natural to give a coloring, which only depends on the norm of the points.

First we choose a suitable prime $p$, two integer parameters $d,l$ and a red/blue coloring of $\mathbb{F}_p$, where $0,d,2d,...,(l-1)d$ 
are colored red and the remaining numbers blue. Then we color each point $x\in\mathbb{E}^n$ to the color of 
$\lfloor |x|^2  \rfloor$ (mod $p$).

Let $X,Y,Z\in\{0,d,2d,...,(l-1)d\}$ be the squared norms of three red points, that form an $\ell_3$, 
already rounded and reduced modulo $p$. Observe that we have $K=2$ here.
To avoid red copies of $\ell_3$, we need two conditions on the parameters $p,d,l$:
\begin{itemize}
    \item $p$ should be large enough so that $p>2(l-1)d+K+2$. This allows us to look at the equations $X-2Y+Z=K'$ for values of $K'$ close to $K$ in the integers instead of modulo $p$, since $-p+K'<X-2Y+Z<p+K'$.
    \item $d\ge 4$ so that $K \leq d-2$. Then $1\le K'\le d-1$. But 
    $K'=X-2Y+Z$ should be $0$ modulo $d$ for a red $\ell_3$, which gives the desired contradiction. 
\end{itemize}

The squared norms of a longer progression $\{x_1,x_2,...,x_n\}$ can be described by a quadratic polynomial function: $|x_i|^2=a i^2+b i+c$. By appropriately rescaling the problem and looking at sub-progressions, we can assume that the coefficients $a$ and $b$ are integral, at the cost of having some error term.  We get $|x_i|^2\approx a(i+b')^2+c'$ and therefore, the function hits every square or every non-square modulo $p$, shifted by a constant. What remains is to choose $l$ large enough, for a fixed $p$, such that $\{C,C+d,C+2d,...,C+(l-1)d\}$ contains both squares and non-squares for every choice of $C$. This guarantees that there are no arbitrarily long blue progressions.

\section{Proof of Theorem \ref{mainthm}}
\label{mainproof}


\begin{lemma}
\label{alpha^2}
Let $x,y,z\in\mathbb{E}^n$ form a configuration congruent to 
$\alpha\ell_3$ 
i.e. $x-2y+z=0$ and $|x-y|^2=|y-z|^2=\alpha^2$. Then we have  $$|x|^2-2|y|^2+|z|^2=2\alpha^2.$$
\end{lemma}

\begin{proof}
 
 \begin{equation*}
\begin{split}   
|x|^2-2|y|^2+|z|^2=& |x|^2-2|y|^2+|2y-x|^2 \\
=&|x|^2-2|y|^2+4|y|^2+|x|^2-4\langle  x , y \rangle \\
=&2|x|^2+2|y|^2-4\langle  x , y \rangle \\ 
=&2|x-y|^2 \\
=&2\alpha^2.
\end{split}
\end{equation*}
\end{proof}

Now we define the red/blue coloring of $\mathbb{E}^n$.
Let $\mathcal{R}$ (resp. $\mathcal{B}$) denote the set of red (resp. blue) points.
Let
$$\mathcal{R}:=\{x\in\mathbb{E}^n\ |\ \left \lfloor |x|^2 \right \rfloor \in\{0,4,8,12\}+29\mathbb{Z} \}$$ and 
$$\mathcal{B}:=\mathbb{E}^n\setminus\mathcal{R}.$$ 

We have to show that there is no red copy of $\ell_3$ and no blue copy of $\ell_{1177}$. Suppose that 
$x,y,z\in\mathbb{E}^n$ form a red configuration congruent to 
$\ell_3$. 

For simplicity, let 
$X=\lfloor |x|^2\rfloor$, 
$Y=\lfloor |y|^2\rfloor$, 
$Z=\lfloor |z|^2\rfloor$ 
and let
$X'=\left\{ |x|^2\right\}$, 
$Y'=\left\{ |y|^2\right\}$, 
$Z'=\left\{ |z|^2\right\}$.

By Lemma \ref{alpha^2}, 
$X+X'-2Y-2Y'+Z+Z'=2$.


\begin{lemma}
We have 
$$X-2Y+Z\in\{1,2,3\}.$$
\end{lemma}

\begin{proof}
For $R=X'-2Y'+Z'$, 
we have $|R|<2$. But then
$X-2Y+Z=2-R$. Since  $X-2Y+Z\in\mathbb{Z}$, $R\in\mathbb{Z}$, therefore, 
$R\in\{-1,0,1\}$, so  $X-2Y+Z\in\{1,2,3\}$.
\end{proof}

Since $x, y, z$ are red, 
$X, Y, Z \in \{0,4,8,12\}+29\mathbb{Z}$. 
To get a contradiction, it is enough 
to show is that the three equations 
$X-2Y+Z=k$, $k\in\{1,2,3\}$,
do not have any solution in $\mathbb{F}_{29}$ such that $X,Y,Z\in\{0,4,8,12\}$. Since $-29+3<-24\leq X-2Y+Z \leq 24 < 29+1$ as inequalities in $\mathbb{Z}$ with  $X,Y,Z\in\{0,4,8,12\}$, it is enough to show that the equations do not have any solutions in $\mathbb{Z}$, which is clear when considered modulo $4$.

\smallskip

Now we show that there is no blue  copy of $\ell_{1177}$.

\begin{lemma}
\label{lemmablue}
For all $c\in\mathbb{F}_{29}$: $$S_{29}+c\not\subseteq \mathbb{F}_{29}\setminus \{0,4,8,12\}.$$
\end{lemma}

\begin{proof}
The squares modulo $29$ are $0$, $1$, $4$, $5$, $6$, $7$, $9$, $13$, $16$, $20$, $22$, $23$, $24$, $25$, $28$, so the non-squares do not contain an arithmetic  
progression of  $4$ elements with consecutive distance $4$, which is equivalent to the statement of the Lemma.
\end{proof}

Suppose that  $\{ x_0,x_1,...,x_{1176}\}$ forms a blue congruent copy of $\ell_{1177}$.
For $0\le i\le 1176$, let $X_i=|x_i|^2$.
By Lemma \ref{alpha^2}, for $0\le i\le 1174$, 
$X_{i+2}=2X_{i+1}-X_i+2$.  We obtain that 
$$X_i=i^2+(X_1-X_0-1)i+X_0, \ \ \  i\in [0,1176].$$ 
Let $\beta:=X_1-X_0-1$. To understand the integral part of $i^2+\beta i+X_0$ we approximate $\beta$ by a rational number.

\begin{lemma}[Dirichlet's approximation theorem \cite{S06}]
For all $\beta\in \mathbb{R}$ and $N\in\mathbb{N}$ there exist $a,d\in\mathbb{Z}$ with $1\leq d\leq N$ such that $$|d\beta-a|\leq \frac{1}{N+1}.$$
\end{lemma}

Let $a,d\in\mathbb{Z}$ satisfy the conditions of Dirichlet's approximation theorem with $N=28$. Then with $\epsilon=\beta-\frac{a}{d}$, 
$|\epsilon|\le\frac{1}{29d}$. So 
$$i^2+\beta i+X_0=i^2+\frac{a}{d}i+X_0+\epsilon i.$$ 
Consider now only every $d$-th 
point, let $X'_j=X_{dj}$. Then 
$$X'_{j}=X_{dj}=(dj)^2+aj+X_0+\epsilon dj,$$ and 
$$|\{j\in\mathbb{N}_0\ |\ dj\leq 1176\}|=1+\bigg\lfloor \frac{1176}{d}\bigg\rfloor\geq 1+\bigg\lfloor \frac{1176}{28}\bigg\rfloor=43.$$

\begin{lemma}
\label{bluethm}
There exists $c\in\mathbb{F}_{29}$ such that 
$$c+S_{29}\subseteq \{\lfloor X'_j \rfloor \;( mod\; 29)\ |\ j\in [0,42]\}.$$
\end{lemma}

\begin{proof}
 \begin{equation*}
\begin{split}   
X'_j=&(dj+(2d)^{-1}a)^2-((2d)^{-1}a)^2+29N+X_0+\epsilon dj \\
=&(dj+s)^2+c'+r+\epsilon dj+29N',
\end{split}
\end{equation*}
where $s,c'\in[0,28]$ and $N,N'$ are suitable integers, $r:=\{ X_0 \}<1$ and $(2d)^{-1}$ is the inverse of $2d$ in $\mathbb{F}_{29}$ considered as an integer in $[0,28]$.

Observe that $dj+s=M$ has a solution $j$ for each $M$ considered as an equation in $\mathbb{F}_{29}$ and that for $\eta:=-sd^{-1}$, $j_-:=\eta-k$ and $j_+:=\eta+k$ with $k\in\mathbb{F}_{29}$, $$(dj_-+s)^2=(dj_++s)^2.$$

Therefore, $$\{(dj+s)^2+c'|j\in[\eta-14,\eta]\}=\{(dj+s)^2+c'|j\in[\eta,\eta+14]\}=c'+S_{29}.$$

Let $\eta_0\in[0,28]$ be a representative of $\eta$.

\begin{itemize}
\item Case 1: $\eta_0<14$

Note that $\eta_0+29\leq 42$. We can assume that either 
$$\lfloor r+\epsilon d\eta_0\rfloor=\lfloor r+\epsilon d(\eta_0+k)\rfloor,\;\;\forall
k\in[0,14],$$ or 
$$\lfloor r+\epsilon d(\eta_0+15)\rfloor=\lfloor r+\epsilon d(\eta_0+k)\rfloor,\;\;\forall
k\in[15,29].$$ 
Indeed, otherwise, as $r+\epsilon dj$ is either increasing or decreasing in $j$, $$|(r+\epsilon d\eta_0)-(r+\epsilon d(\eta_0+29))|>1,$$ which contradicts $|29\epsilon d|\leq 1$.

But then either 
$$\{\lfloor X'_j \rfloor \;( mod\; 29)\ |\ j\in [\eta_0,\eta_0+14]\}=c'+\lfloor r+\epsilon d\eta_0\rfloor+S_{29}$$ or 
$$\{\lfloor X'_j \rfloor \;( mod\; 29)\ |\ j\in [\eta_0+15,\eta_0+29]\}=c'+\lfloor r+\epsilon d(\eta_0+15)\rfloor+S_{29}.$$

\item Case 2: $\eta_0\geq 14$

We can proceed similarly.  Note that $\eta_0-14\geq 0$, $\eta_0+14\leq 42$.
We can assume that either 
$$\lfloor r+\epsilon d\eta_0\rfloor=\lfloor r+\epsilon d(\eta_0+k)\rfloor,\;\;\forall
k\in[-14,0]$$ or 
$$\lfloor r+\epsilon d\eta_0\rfloor=\lfloor r+\epsilon d(\eta_0+k)\rfloor,\;\;\forall
k\in[0,14],$$ 
otherwise, as $r+\epsilon dj$ is either increasing or decreasing in $j$, $$|(r+\epsilon d(\eta_0-14))-(r+\epsilon d(\eta_0+14))|>1,$$ which contradicts $|29\epsilon d|\leq 1$.

But then either 
$$\{\lfloor X'_j \rfloor \;( mod\; 29)\ |\ j\in [\eta_0-14,\eta_0]\}=c'+\lfloor r+\epsilon d\eta_0\rfloor+S_{29}$$ or 
$$\{\lfloor x'_j \rfloor \;( mod\; 29)\ |\ j\in [\eta_0,\eta_0+14]\}=c'+\lfloor r+\epsilon d\eta_0\rfloor+S_{29}.$$

\end{itemize}

\end{proof}

Lemmas \ref{lemmablue} and \ref{bluethm} together imply that there is no blue copy of $\ell_{1177}$ in the construction and that completes the proof of Theorem \ref{mainthm}.
$\Box$

\section{Proof of Theorem \ref{generalthm}}
\label{generalproof}

Analogously to the distance one case, set 
$$\mathcal{R}:=\{x\in\mathbb{E}^n\ |\ \lfloor |x|^2 \rfloor \in\{0,5,10,15,20\}+47\mathbb{Z} \}$$ and $$\mathcal{B}:=\mathbb{E}^n\setminus\mathcal{R}.$$ 

We will study for which $\alpha_{red}\in\mathbb{R}$, is it true that
there are no three red points $x, y, z\in \mathbb{E}^n$ that satisfy the equation $|x|^2-2|y|^2+|z|^2=2\alpha_{red}^2$ and for which $\alpha_{blue}\in\mathbb{R}$ is it true that
there are no $6628$ blue points $x_0,x_1,...,x_{6627}$ in $\mathbb{E}^n$ that satisfy the equations $|x_i|^2-2|x_{i+1}|^2+|x_{i+2}|^2=2\alpha_{blue}^2$ for $i\in [0,6625]$ and then give the corresponding ratios $\alpha_{blue}/\alpha_{red}$.
We choose the the bigger prime $47$ and the above coloring as it allows us some freedom to increase the distances in the $3$-term arithmetic progression.

\begin{lemma}
\label{lemmascale}
Let $x, y, z$ form  a copy of $\alpha_{red}{\ell}_3$ with $47N+1\leq\alpha_{red}^2\leq 47N+3/2$ and $N\in\mathbb{Z}_{\geq 0}$ then 
$$\lfloor |x|^2\rfloor-2\lfloor |y|^2\rfloor+\lfloor |z|^2 \rfloor\in \{1,2,3,4\} \;( mod\; 47).$$
\end{lemma}

\begin{proof}
Let $R:=(\{|x|^2\}-2\{|y|^2\}+\{|z|^2\})$, clearly $|R|<2$. Now

$\lfloor|x|^2\rfloor-2\lfloor |y|^2\rfloor+\lfloor |z|^2\rfloor=2\alpha^2-R.$ 
Since $2\alpha^2\in [2,3]+47\mathbb{Z}$ and  
$\lfloor |x|^2\rfloor-2\lfloor |y|^2\rfloor+\lfloor |z|^2\rfloor\in\mathbb{Z}$, 
we have $\lfloor|x|^2\rfloor-2\lfloor |y|^2\rfloor+\lfloor |z|^2\rfloor\in \{1,2,3,4\} \;( mod\; 47).$

\end{proof}

\begin{lemma}
For all $c\in\mathbb{F}_{47}$: $$S_{47}+c\not\subseteq \mathbb{F}_{47}\setminus \{0,5,10,15,20\},$$ and 
$$(\mathbb{F}_{47}\setminus S_{47}^*)+c\not\subseteq \mathbb{F}_{47}\setminus \{0,5,10,15,20\},$$

\end{lemma}

\begin{proof}
Assume for a contradiction that 
$L:=\{c, c+5, c+10, c+15, c+20\}$ is either contained in  
$(\mathbb{F}_{47}\setminus S_{47})=(\mathbb{F}^*_{47}\setminus S_{47}^*)$ or in  $S_{47}^*$. 
The squares in $\mathbb{F}_{47}^*$ are $1$, $2$, $3$, $4$, $6$, $7$, $8$, $9$, $12$, $14$, $16$, $17$, $18$, $21$, $24$, $25$, $27$, $28$, $32$, $34$, $36$, $37$, $42$ and it is easy to see that there is no interval and no gap of size $5$ in the squares. Now $$\{19c,19c+1,19c+2,19c+3,19c+4\}=19L\subseteq19(\mathbb{F}^*_{47}\setminus S_{47}^*)=S_{47}^*$$ or $$19L\subseteq 19S_{47}^*=(\mathbb{F}^*_{47}\setminus S_{47}^*),$$ a contradiction.
\end{proof}

We now choose a longer arithmetic progression, compared to the distance one case which allows us to add an additional error term.

Let $\alpha_{blue}\in\mathbb{R}_+$ such that $\alpha_{blue}^2=b+\epsilon_2$ with 
$b\in\mathbb{N}\setminus 47\mathbb{Z}$ and 
$0<\epsilon_2 <17(7\cdot 47^4 \cdot 48)$ and suppose that 
$\{ x_0,x_1,...,x_{8648}\}$ forms a copy of $\alpha_{blue}{\ell}_{8649}$.
For $0\le i\le 8648$, let $X_i=|x_i|^2$. By lemma \ref{alpha^2}, for 
$0\le i\le 8646$,
$X_{i+2}=2X_{i+1}-X_i+2\alpha_{blue}^2$ consequently $X_i=\alpha_{blue}^2i^2+(X_1-X_0-1)i+X_0$. 
Let $\beta:=X_1-X_0-1$.  
Let $a,d\in\mathbb{Z}$ satisfy the conditions of Dirichlet's approximation theorem with $N=47$, so $$\alpha_{blue}^2i^2+\beta i+x=\alpha_{blue}^2i^2+\frac{a}{d}i+x+\epsilon_1 i$$ with $|\epsilon_1|\leq 1/(48d)$. Consider only every $d$-th point, let $X'_j=X_{dj}$. Then 
$$X'_{j}:=X_{dj}=\alpha_{blue}^2(dj)^2+aj+X_0+\epsilon_1dj,$$ 

and $$|\{j\in\mathbb{N}_0\ |\ dj\leq 8648\}|=1+\bigg\lfloor \frac{8648}{d}\bigg\rfloor\geq 1+\bigg\lfloor \frac{8648}{47}\bigg\rfloor=185.$$

\begin{lemma}
\label{bluethmgeneral}
There exists $c\in\mathbb{F}_{47}$ such that $$c+S_{47}\subseteq \{\lfloor X'_j \rfloor \;( mod\; 47)\ |\ j\in [0,184]\}$$ or 
$$c+(\mathbb{F}_{47}\setminus S_{47}^*)\subseteq \{\lfloor X'_j \rfloor \;( mod\; 47)\ |\ j\in [0,184]\}$$
\end{lemma}

\begin{proof}

Suppose first that $d<47$. Then

\begin{equation*}
\begin{split}   
X'_j=&\alpha_{blue}^2d^2j^2+aj+X_0+\epsilon_1dj \\
=&bd^2j^2+aj+X_0+\epsilon_1dj+\epsilon_2d^2j^2 \\
=&b(dj+(2db)^{-1}a)^2-((2b)^{-1}a)^2+47N+X_0+\epsilon_1dj+\epsilon_2d^2j^2 \\
=&b(dj+s)^2+c'+r+\epsilon_1dj+\epsilon_2d^2j^2+47N',
\end{split}
\end{equation*}
where $s,c'\in[0,46]$ and $N,N'$ are suitable integers, $r:=\{ X_0 \}<1$ and $(2db)^{-1}$ is the inverse of $2db$ in $\mathbb{F}_{47}$ considered as an integer in $[0,46]$.

Now either $b$ is a square and $b(dj+s)^2$ will run through $S_{47}$ or otherwise
$b(dj+s)^2$ will run through $(\mathbb{F}_{47}\setminus S_{47}^*)$. In both cases denote the corresponding set by $S$.

Observe that $dj+s=M$ has a solution $j$ for each $M$ considered as an equation in $\mathbb{F}_{47}$ and that for $\eta:=-sd^{-1}$, $j_-:=\eta-k$ and $j_+:=\eta+k$ with $k\in\mathbb{F}_{47}$, $$b(dj_-+s)^2=b(dj_++s)^2.$$

Therefore, $$\{b(dj+s)^2+c'\ |\ j\in[\eta-23,\eta]\}=\{b(dj+s)^2+c'\ |\ j\in[\eta,\eta+23]\}=c'+S.$$

Let  $\eta_0$ be a representative of $\eta$ in $[0,46]$ and write $E(\mu):=r+\epsilon_1d\mu+\epsilon_2d^2\mu^2$. As $E$ is at most quadratic, it changes its monotonicity at most once. 

\begin{itemize}
\item Case 1: $\eta_0<23$.

Let $\delta:=0$ if $E$ is monotone in $[\eta_0,\eta_0+47]$ and $\delta:=47$ otherwise. 
Let $\eta_1:=\eta_0+\delta$. Now $E$ is monotone in $I:=[\eta_1,\eta_1+47]$.

Note that $\eta_1+47\leq \delta+69\leq 116<185$. Assume that $$\lfloor E(\eta_1)\rfloor=\lfloor E(\eta_1+k)\rfloor,\;\;\forall
k\in[0,23]$$ or $$\lfloor E(\eta_1+24)\rfloor=\lfloor E(\eta_1+k)\rfloor,\;\;\forall
k\in[24,47].$$ 

Then $$\{\lfloor X'_j \rfloor \;( mod\; 47)\ |\ j\in [\eta_1,\eta_1+23]\}=c'+\lfloor E(\eta_1)\rfloor+S$$ or $$\{\lfloor X'_j \rfloor \;( mod\; 47)\ |\ j\in [\eta_1+24,\eta_1+47]\}=c'+\lfloor E(\eta_1+24)\rfloor+S.$$

Otherwise, as $E$ is monotone in $I$, $$|E(\eta_1)-E(\eta_1+47)|>1,$$ which is a contradiction to 
\begin{equation*}
\begin{split}   
&|E(\eta_1)-E(\eta_1+47)| \\
=&|47\epsilon_1d+(94\eta_1+47^2)\epsilon_2d^2| \\
\leq&|47\epsilon_1d|+|(94\eta_1+47^2)\epsilon_2d^2| \\
\leq&\frac{47}{48}+(94\cdot 93+47^2)46^2|\epsilon_2| \\
\leq&\frac{47}{48}+5\cdot 47^4|\epsilon_2|\leq 1 .
\end{split}
\end{equation*}

\item Case 2: $\eta_0\geq 23$

Similarly, let $\delta:=0$ if $E$ is monotone in $[\eta_0-23,\eta_0+23]$ and $\delta:=47$ otherwise. 
Let $\eta_1:=\eta_0+\delta$. Now $E$ is monotone in $I:=[\eta_1-23,\eta_1+23]$.

Note that $\eta_1-23\geq \delta\geq 0$, $\eta_1+23\leq \delta+69\leq 116<185$ and assume that $$\lfloor E(\eta_1)\rfloor=\lfloor E(\eta_1+k)\rfloor,\;\;\forall
k\in[-23,0]$$ or $$\lfloor E(\eta_1)\rfloor=\lfloor E(\eta_1+k)\rfloor,\;\;\forall
k\in[0,23].$$ 

Then $$\{\lfloor x'_j \rfloor \;( mod\; 47)\ |\ j\in [\eta_1-23,\eta_1]\}=c'+\lfloor E(\eta_1)\rfloor+S$$ or $$\{\lfloor x'_j \rfloor \;( mod\; 47)\ |\ j\in [\eta_1,\eta_1+23]\}=c'+\lfloor E(\eta_1)\rfloor+S.$$

Otherwise, as $E$ is monotone in $I$, $$|E(\eta_1-23)-E(\eta_1+23)|>1,$$ which is a contradiction to 
\begin{equation*}
\begin{split}   
&|E(\eta_1-23)-E(\eta_1+23)| \\
=&|46\epsilon_1d+(92\eta_1+2\cdot23^2)\epsilon_2d^2| \\
\leq&|46\epsilon_1d|+|(92\eta_1+2\cdot23^2)\epsilon_2d^2| \\
\leq&\frac{46}{48}+(92\cdot 93+2\cdot23^2)46^2|\epsilon_2| \\
\leq&\frac{46}{48}+5\cdot 47^4|\epsilon_2|< 1 .
\end{split}
\end{equation*}

\end{itemize} 

Now let $d=47$, then $a\neq 0$. Then 

\begin{equation*}
\begin{split}   
X'_j=&\alpha_{blue}^2d^2j^2+aj+X_0+\epsilon_1dj \\
=&aj+X_0+\epsilon_1dj+\epsilon_2d^2j^2+47N \\
=&aj+c+r+\epsilon_1dj+\epsilon_2d^2j^2+47N',
\end{split}
\end{equation*}

where $c\in [0,46]$, $N,N'$ are integers and $r:=\{ X_0 \}<1$.

Again, write $E(\mu):=r+\epsilon_1d\mu+\epsilon_2d^2\mu^2$ and let $\delta:=0$ if $E$ is monotone in $[0,92]$ or $\delta:=92$ otherwise. Now $E$ is monotone in $I:=\delta+[0,92]$ and $I\subseteq [0,184]$.
If $$\lfloor E(\delta)\rfloor=\lfloor E(\delta+k)\rfloor,\;\;\forall
k\in[0,46]$$ 
then 
$$\{\lfloor X'_j \rfloor \;( mod\; 47)\ |\ j\in [\delta,\delta+46]\}=\mathbb{F}_{47}\supseteq S_{47}.$$ 
Otherwise, let $\eta\in \delta+[1,46]$ be such that $\lfloor E(\eta-1)\rfloor=\lfloor E(\eta)\rfloor$. Note that $[\eta-1,\eta+46]\subseteq I$. 
Again, if $$\lfloor E(\eta)\rfloor=\lfloor E(\eta+k)\rfloor,\;\;\forall
k\in[0,46]$$ 
then 
$$\{\lfloor X'_j \rfloor \;( mod\; 47)\ |\ j\in [\eta,\eta+46]\}=\mathbb{F}_{47}\supseteq S_{47}.$$ Otherwise, as $E$ is monotone in $I$, $$|E(\eta-1)-E(\eta+46)|>1,$$ which is a contradiction to 
\begin{equation*}
\begin{split}   
&|E(\eta-1)-E(\eta+46)| \\
=&|47\epsilon_1d+(94\eta+46^2-1)\epsilon_2d^2| \\
\leq&|47\epsilon_1d|+|(94\eta+46^2-1)\epsilon_2d^2| \\
\leq&\frac{47}{48}+(94\cdot 140+47^2)46^2|\epsilon_2| \\
\leq&\frac{47}{48}+7\cdot 47^4|\epsilon_2|\leq 1 .
\end{split}
\end{equation*}

\end{proof}


\subsection{Large numbers}

Together, Lemmas \ref{lemmascale} and \ref{bluethmgeneral} prove Theorem \ref{generalthm} for $\alpha^2\geq 2$. The intervals $[2n/3,n]$ for $n\in\mathbb{Z}_{\geq 3}\setminus 47\mathbb{Z}$ cover $\mathbb{R}_{\geq 2}$ as $2n/3\leq n-1$ for $n\geq 3 $ and $2n/3\leq n-2$ for $n\in 1+ 47\mathbb{N}$. Therefore it is sufficient to use red progressions 
$\alpha_{red}\ell_3$ and blue progressions $\alpha_{blue}\ell_m$ with $1\leq\alpha_{red}^2\leq 3/2$ and $\alpha_{blue}^2\in\mathbb{Z}_{\geq 3}\setminus 47\mathbb{Z}$.

\subsection{Small numbers}

Let $c:=7\cdot47^4\cdot48$. Analogously to the case of large ratios, Lemmas \ref{lemmascale} and \ref{bluethmgeneral} prove Theorem \ref{generalthm} for $\alpha^2\leq 1/(2c)$. The intervals $[n/(1+1/c),n]$ for $n\in\mathbb{Z}_{\geq 2(c+1)}\setminus 47\mathbb{Z}$ cover $\mathbb{R}_{\geq 2c}$ as $n/(1+1/c)\leq n-2$ for $n\geq 2(c+1)$. Therefore it is sufficient to use red progressions 
$\alpha_{red}\ell_3$ and blue progressions 
$\alpha_{blue}\ell_m$ with $1\leq\alpha_{blue}^2\leq 1+1/c$ and $\alpha_{red}^2\in\mathbb{Z}_{\geq 2(c+1)}\setminus 47\mathbb{Z}$.

\subsection{Rational numbers}
\label{rationalsubsection}

\begin{itemize}
\item Case 1:

Let $\alpha^2=p/q$ where $p,q\in\mathbb{N}$, $47\nmid pq$ and $k,r\in\mathbb{Z}$ such that $0<r<47$ and $q=47k+r.$ Let $r^{-1}$ be the inverse of $r$ in $\mathbb{F}_{47}$ considered as an integer in $[1,46]$ and let $\alpha_{red}:=\sqrt{qr^{-1}}$ and $\alpha_{blue}:=\sqrt{pr^{-1}}$. Then our construction does not contain any red copy of $\alpha_{red}\ell_3$ and no blue copy of $\alpha_{blue}\ell_m$ which proves Theorem \ref{generalthm} for $\alpha=\alpha_{blue}/\alpha_{red}$.

\item Case 2:

Let $\alpha^2=p/q$ where $p,q\in\mathbb{N}$, $47\mid p$ and $\gcd(p,q)=1$. Let $g$ be the inverse of $q$ in $\mathbb{Z}/47p\mathbb{Z}$ considered as an integer in $[0,47p-1]$ and $N\in\mathbb{Z}_{\geq0}$ such that $gq=47pN+1$. Then $$(p+1)g\frac{q}{p}=\frac{(p+1)(47pN+1)}{p}=1+\frac{1}{p}+47N(p+1).$$ Now, since $47\nmid (p+1)g$ and $0\leq 1/p\leq 1/2$, our construction does not contain any red copy of $\alpha_{red}\ell_3$ and no blue copy of $\alpha_{blue}\ell_m$ for $\alpha_{red}:=\sqrt{(p+1)gq/p}$ and $\alpha_{blue}:=\sqrt{(p+1)g}$, which proves Theorem \ref{generalthm} for $\alpha=\alpha_{blue}/\alpha_{red}$.

\end{itemize}

\subsection{Irrationals}

Let $\alpha^2$ be irrational. Then the sequence $((47k+1)/(47\alpha^2))_{k\in\mathbb{N}}\;( mod\; 1)$ is uniformly distributed (see e.g. \cite[Theorem~3.2]{KuipersNierderreiter}) and  equivalently $((47k+1)/\alpha^2)_{k\in\mathbb{N}}\;( mod\; 47)$ is uniformly distributed. In particular there exist $p,q\in\mathbb{N}$ with $p\equiv q\equiv 1 \;( mod\; 47)$ and  such that $$q<\frac{p}{\alpha^2}<q+\frac{1}{2}.$$

Now let $\alpha_{red}:=\sqrt{p/\alpha^2}$ and $\alpha_{blue}:=\sqrt{p}$. Then our construction does not contain any red copy of $\alpha_{red}\ell_3$ and no blue copy of $\alpha_{blue}\ell_m$ which proves Theorem \ref{generalthm} for $\alpha=\alpha_{blue}/\alpha_{red}$.

\subsection*{Remarks}
We believe that the bound $1177$ in Theorem \ref{mainthm} is far from optimal.
Let $m_3$ be the largest integer for which  
$\mathbb{E}^n\rightarrow (\ell_3,\ell_{m_3})$ holds for every $n\ge 2$, or equivalently,
$\mathbb{E}^2\rightarrow (\ell_3,\ell_{m_3})$ holds.
Erd\H os et al. \cite{EGM2} showed that
$\mathbb{E}^2\rightarrow (\ell_4,\ell_2)$, and recently Currier et al. \cite{cmh2024-3prog} 
showed that
$\mathbb{E}^2\rightarrow (\ell_3,\ell_3)$,
therefore, 
$m_3\ge 3$. 
We can not rule out the possibility that $m_3=3$. 

Similarly, let $m_4$ (resp. $m_5$) be the largest integer for which  
$\mathbb{E}^n\rightarrow (\ell_4,\ell_{m_4})$ 
(resp. $\mathbb{E}^n\rightarrow (\ell_5,\ell_{m_5})$)
holds for every $n\ge 2$. Clearly, $m_5\le m_4\le m_3$. 
Improving the above result of Erd\H os et al., 
Tsaturian \cite{T17} showed that 
$\mathbb{E}^2\rightarrow (\ell_5,\ell_2)$, consequently, $m_4\ge m_5\ge 2$. 

For $\ell_6$, the situation is different. Define $m_6$ analogously. It is not known, whether
$m_6\ge 2$, on the other hand, Erd\H os et al. \cite{EGM1} proved that
$\mathbb{E}^n\nrightarrow (\ell_6,\ell_6)$ for every $n\ge 2$, therefore, $m_6\le 5$.

\smallskip

For Theorem \ref{generalthm} 
we also believe that the bound 
$8649$ is far from optimal and the conditions for $\alpha$ can be dropped. It is clear that by using different primes in the proof we get a finite bound on $m$ for every value of $\alpha$:
The P\'olya–Vinogradov inequality  (see e.g. \cite[p.~135]{Davenport}) guarantees that the length of gaps in the squares is sub-linear in $p$. We can therefore use the construction in subsection~\ref{rationalsubsection} with any high enough prime that does not divide $q$. This will however not give a uniform bound for $m$.

\smallskip

During the review progress of this manuscript Currier, Moore and Yip published a preprint~\cite{currier2024avoidingshortprogressionseuclidean} on the same problem. In particular it contains an improvement to Theorem~\ref{mainthm}.

\subsection*{Acknowledgements}
 J. F. was supported by the Austrian Science Fund (FWF) under
the project W1230. G. T. was supported by the National Research, Development
and Innovation Office, NKFIH, K-131529 and ERC Advanced Grant “GeoScape,”
No. 882971.
The authors also thank the anonymous referee for comments on this manuscript.

\printbibliography

\end{document}